# Additive Separability, Optimization, and Trivial Webs

Thomas Russell

Santa Clara University

Dept. of Economics

The author wished to thank Professor J B Cooper for many helpful discussions on this topic. He also wishes to thank Professor A Kushner for helpful comments. All errors are his own.

**Introduction:**

It has long been understood that there is an intimate connection between the additively separable functions of economics and trivial 3-webs on the plane, Debreu (1959), (1960). Recalling that on the plane a 3-web is given by three foliations of smooth curves in general position, two such webs are said to be locally equivalent at a point p if there exists a local diffeomorphism on a neighborhood of p which exchanges them.
If we consider a 3-web on the *x,y* plane given by

a)  a family of indifference curves $f(x, y) = const$
b)  the horizontals $y = const$ and  (1)
c)  the verticals $x = const$

then it is clear that if we can write the function $f$ in the form

$$\varphi(f()) = U_1(x) + U_2(y)$$

the web is locally equivalent to the trivial web defined as

$$X = const, \; Y = const, \; X + Y = const.$$

This fact forms the heart of Debreu's use of the Blaschke-Bol-Thomsen "hexagon closing" condition for triviality to characterize additively separable preferences, see the more extensive treatments in Wakker (1989) and Vind (2002).

In this paper we show that web triviality provides a geometric test for the prior and deeper "integrability" question of whether or not an observed family of price reaction functions is generated by a maximization process, Samuelson (1972). The paper draws on a theorem of Tabachnikov (1993) which characterizes local triviality equivalence for both 3-webs and (under area preserving diffeomorphisms called symplectomorphisms) for 2-webs. Tabachnikov's result is remarkable (at least to economists) in that it provides in one theorem a geometric condition for both additive separability and integrability. For the latter problem we note that a 3-web on the plane, given by the verticals, the horizontals, and the level curves of a certain transformation function is geometrically trivial. This result is related to a fact first noted by Samuelson (1960), see also Samuelson (1972), Samuelson (1983), Cooper, Russell, and Samuelson (2001), and Cooper and Russsell (2006)  that the additive separability of  a certain  Jacobian  guarantees integrability.

**Trivial 3-webs: Additive Separability and the Chern Connection.** As Tabachnikov (op.cit.)  notes, an elegant characterization of trivial webs can be given in terms of the curvature of a connection associated to it: for example, a 3-web on the plane is trivial if and only if the curvature of the Chern connection vanishes. For the 3-web given in (1), K, the curvature of the Chern connection is given by

$$K = -\frac{1}{f_x f_y} (\log(\frac{f_x}{f_y}))_{xy} \qquad (2)$$

Akivis et al  (2005). Here subscripts denote partial derivatives.
In the case of decision making under uncertainty, it has already been noted, Russell (2003), that Chern curvature measures the extent to which individuals fail to satisfy the axioms of expected utility maximization. In particular, when K=0,  i.e. the web is trivial, the right hand side of (2) set equal to zero gives the de Saint Robert equation (Vind (op.cit.) p 79), and when the level curves of an individual over mixtures of sure things satisfy this equation, the individual's preferences can be represented by a function $V = \frac{1}{2}((U(x)+U(y))$. Further discussion of this equation can be found in Cooper Russell and Samuelson (2004).

**Optimizing Behavior and Trivial Webs.** The concept of a connection associated to a web can be extended to other webs. Once a connection is defined, its curvature can be calculated. We now show how zero curvature of a canonical connection associated with a 2-web implies integrability.

We follow Samuelson (1972) in considering a monopolist facing a downward sloping demand function for output who hires two inputs $q_1, q_2$ at competitive prices $p_1, p_2$.
In Fig. 1 we show the factor demands of this firm in two regimes of constraint. The flatter (bold) family of curves shows $q_1=f_1(p_1; p_2)$, the demand for factor 1 as a function of its own price, holding fixed the price of factor 2. A higher price $p_2$ moves the demand function up and to the right. The second (steeper) family of curves shows $q_1=f_2(p_1; q_2)$, the demand for factor 1 as a function of its own price, holding fixed the quantity of factor 2. Giving the firm a larger allocation of $q_2$ moves the demand function for $q_1$ down to the left. The integrability question is whether or not this firm is profit maximizing.

$p_1$

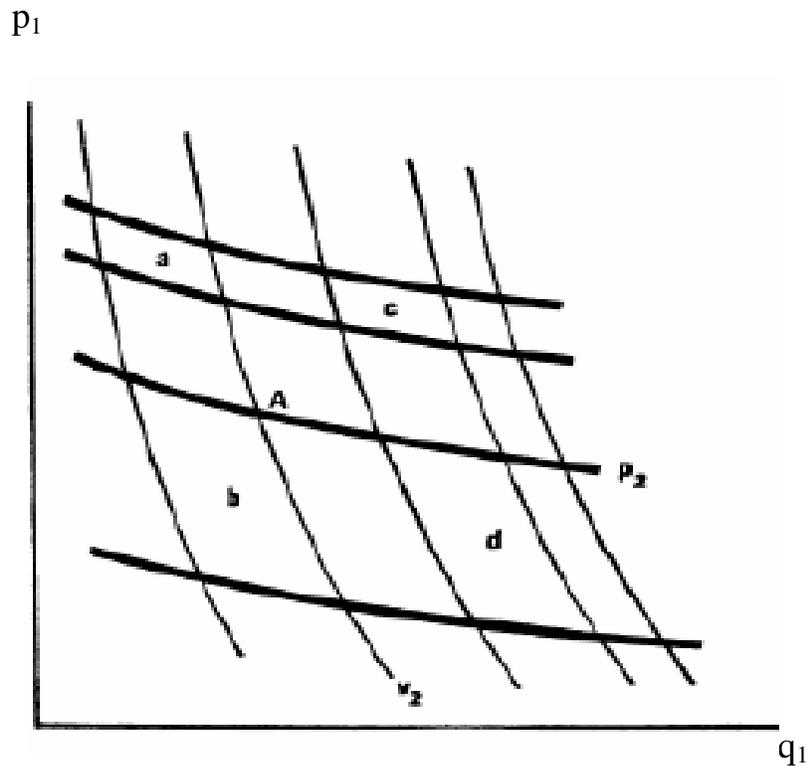

$q_1$

**Figure 1 : The Constrained Demand Function Web**

Samuelson (1972) does not use the language of web geometry. Indeed his discussion predates the use of differential geometry to describe the equilibrium thermodynamics to which his diagram is mathematically equivalent[1]. Therefore before discussing Tabachnikov's theorem, it may be useful to reframe the economic problem in the modern language of contact and symplectic geometry. [2]

Writing $\Pi$ for profits, the data of this problem can be viewed as a point in $R^5$ with independent coordinates $\Pi, q_1, q_2, p_1, p_2$. By standard duality theory, if the firm is profit maximizing, there will exist a profit function relating these variables through a differential form

$$d\Pi = -q_1 dp_1 - q_2 dp_2$$

That is to say, if we introduce the form

$$\omega = d\Pi + q_1 dp_1 + q_2 dp_2 \qquad (2)$$

the set of $q,p$ combinations consistent with a profit function are the zeros of $\omega$.

An equation of the form $\omega = 0$ is a Pfaffian equation. Its solution is a sub-manifold of $R^5$ whose dimension is determined by an examination of the exterior derivatives of $\omega$. Calculating the appropriate terms, we have

$$d\omega = dq_1 \wedge dp_1 + dq_2 \wedge d_2 p$$
$$d\omega \wedge d\omega = 2 dq_1 \wedge dp_1 \wedge dq_2 \wedge dp_2$$
$$d\omega \wedge d\omega \wedge d\omega = 2 d\Pi \wedge dq_1 \wedge dp_1 \wedge dq_2 \wedge dp_2 \neq 0$$

Thus $(R^5, \omega)$, is a contact vector space i.e. a manifold on which there is defined a form $\omega$ for which $d\omega \wedge (d\omega)^2 \neq 0$. The form $\omega$ generates a 4-dimensional non-integrable distribution, the so-called **contact distribution.** Two-dimensional integral manifolds of the distribution are called *Legendrian* **manifolds.** If we project out the $\Pi$ dimension, we obtain an even dimensional manifold on which there is defined a symplectic form $d\omega = dq_1 \wedge dp_1 + dq_2 \wedge d_2 p$. Sub-manifolds of symplectic manifolds for which $d\omega = 0$ are known as **Lagrangian** sub-manifolds. Thus the optimizing processes of economics form a Legendrian (resp. Lagrangian) sub-manifolds of contact (resp. symplectic) manifolds given by the original data[3].

---

[1] In thermodynamics this classic diagram with appropriate physical coordinates was first used by Watt in the early 19th century to describe steam engines. I is known as the Watt indicator diagram.

[2] The modern treatment of thermodynamics goes back to Hermann (1973) and Arnold (1990). In the last 20 years there have been a number of significant contributions form Mrugala and his co-authors, e.g. Mrugala (1991) (2000)

Choosing the Lagrangian formulation, the question thus becomes one of identifying Lagrangian submanifolds in terms of the restricted input demand functions shown in Fig. 1.

**Identifying Lagrangian Sub-manifolds.**

As Samuelson noted, if we consider general "quadrilaterals" cut off by the "threads" of the two foliations, for example in Fig. 1 areas such as those marked a,b,c,d, the equality ab = cd, provides a necessary condition for profit maximization.. Call this the Samuelson (S) area ratio condition. In what way is this related to Lagrangian sub-manifolds.? Associated with any $2n$ –dimensional symplectic vector space is the $n(n+1)/2$ –dimensional manifold of its Lagrangian subspaces, the so called **Lagrangian Grassmannian.** We may use the data provided by the input demand functions to generate one candidate Lagrangian submanifold, the graph of a function over $p_1, q_1$. When this surface is a Lagrangian submanifold, ( i.e. $d\omega$ is to equal zero on S) the associated map from the $p_1, q_1$-plane to the $p_2, q_2$-plane must clearly be area preserving but orientation reversing, (since $d\omega = 0$ implies $dq_1 \wedge dp_1 = -dq_2 \wedge dp_2$, ). In terms of web geometry, then, the question is whether or not the 2-web in Fig 1 is equivalent to the trivial 2-web $q_2 = constant$  $p_2 = constant$ under an area preserving map. Orientation reversal is guaranteed by the fact that the leaves of the $q_2$ foliation are already inverse ordered.

But it is known that

---

[3] Obviously much more may be said regarding contact geometry and economics. In Russell (2007) some of this structure is presented. Lichnerowicz (1970) provides the first use of these techniques in economics of which we are aware. Salamon et al (2006) provides a clear overview of the role of contact and symplectic geometry in thermodynamics. Making the association: conjugate variables in thermodynamics = dual variables in economics, this paper may be considered a template for economics

Result (Tabachnikov) When the S condition is satisfied, input demand functions are locally equivalent to the trivial web $q_2=constant$ $p_2=constant$ under an area preserving (and in this case) orientation reversing map.

It follows immediately that

Corollary: The surface generated in ($R^4 d\omega$) by the graph over the $p_1, q_1$-plane associated with the restricted input demand functions is a Lagrangian sub-manifold.

**Profit Maximization and the Hess Connection**: As we noted earlier, we can associate with a 3-web a canonical connection (the Chern connection) whose flatness guarantees additive separability. In this case also we can attach to the 2-web of input demand functions a canonical connection whose flatness is a necessary condition for profit maximization.

This connection was introduced by Hess (1978) as a canonical connection of bi Lagrangian manifolds. To see its relevance to maximizing processes, we need to introduce some further terminology.

Definition: A *Lagrangian foliation* is a foliation of a symplectic manifold on whose leaves the symplectic form vanishes.

Definition A *bi-Lagrangian manifold* is a symplectic
manifold endowed with two transversal Lagrangian foliations)

From these definitions it follows that any pair of transverse foliations of the 2-dimensional $p_1, q_1$-plane with the standard symplectic form is a bi-Lagrangian manifold, so in particular we have

Result: The restricted input demand functions of Fig.1 constitute a bi-Lagrangian manifold.

Finally we introduce the Hess connection $\nabla$ on this manifold

Definition: The *Hess connection* of a bi-Lagrangian manifold is the unique symplectic connection $\nabla$ which parallelizes the two foliations $F1$ and $F2$.

We have the following theorem

Theorem (Hess) If the curvature of $\nabla$ vanishes identically on an n-dimensional manifold, , then the bi_Lagrangian manifold is locally isomorphic to $R^n$ with the standard symplectic structure.

As Tabachnikov notes, the failure of the Samuelson area condition is the obstruction to the vanishing of the Hess connection on the plane. Thus when the Samuelson area ratio condition is satisfied, the curvature of the Hess connection vanishes and there is an area preserving local diffeomorphism which takes the 2-web generated by the two input demand foliations into the trivial 2-web. In other words, the vanishing of the Hess connection is a necessary condition for profit maximization.

**Additive Separabilty and Integrability.** The problem of additive separability discussed by Debreu (op. cit.) and the problem of integrability discussed by Samuelson (op. cit.) are conceptually distinct. It seems little more than a coincidence that to each problem we can associate a (different) web whose canonical connection is flat when the relevant economic hypotheses are satisfied.

There is, however, a mathematical link between the two problems which goes beyond mere happenstance. From the work of Libermann and Marle (1987) p 140, we know that every point in Fig 1 has a neighborhood on which $q_2, p_2$ may be used as local coordinates and for which the local expression for $\omega$ is given by

$$\omega = a(q_1, p_1) dq_1 \wedge dp_1 \qquad (3)$$

Moreover there exists a differentiable function $S$ which satisfies

$$a(q_1, p_1) = \frac{\partial^2 S}{\partial dq_1 dp_1} \qquad (4).$$

Integrability requires that $dq_1 \wedge dp_1 = dq_2 \wedge dp_2$.

The area of a rectangle[4] with vertices (0,0), (ε,0), (ε, δ),(0, δ) is

A():=S (ε, δ)-S(ε,0) -S(0, δ) +S(0,0).

Using the Samuelson area ratio condition equality

---

[4] The following argument is due to Tabachnikov, private communication.

$$A(\varepsilon, \delta) A(-\varepsilon,- \delta) = A(\varepsilon,- \delta) A(-\varepsilon, \delta),$$

the Taylor expansion in $\varepsilon, \delta$ to 4-th order gives

$$S_{\{q_1,p_1\}} S_{\{q_1,q_1, p_1 p_1\}} = S_{\{q_1,q_1,p_1\}} S_{\{q_1,p_1 p_1\}}.$$

This equation is easily solved: $(\ln S_{\{q_1 p_1\}})_{\{q_1 p_1\}}=0$, and hence $S_{\{q_1 p_1\}}=f(q_1)g(p_1)$. Thus our area form is $f(q_1)g(p_1)\, dq_1 dp_1$ and the diffeomorphism $q_2 = F(q_1)$, $p_2 = G(p_1)$ where F'=f, G'=g gives $dq_1 \wedge dp_1 = dq_2 \wedge dp_2$. The Samuelson area condition guarantees that

$$\ln a(q_1, p_1) = \ln \frac{\partial^2 S}{\partial dq_1 dp_1} \text{ splits as}$$

$$\ln a(q_1, p_1) = \ln \frac{\partial^2 S}{\partial dq_1 dp_1} = f(q_1) + g(p_1).$$

It therefore can be seen that the hexagon condition applied to a 3-web given by the level curves of a certain function associated with the symplectic 2-web, the verticals and the horizontals is a necessary condition for the integrability of a family of input demand functions. A similar argument noting the splitting of the Jacobian of the transformation was given by Samuelson (1960), see also Cooper Russell Samuelson(2001).

**Conclusion :** In this paper we have shown that two seemingly unrelated problems in economics, the hypothesis of integrability and the hypothesis of additive separability are linked by the absence of curvature of connections on webs naturally associated with each problem.